\documentclass[12pt]{article}
\usepackage{amssymb}
\usepackage{latexsym,bm}
\usepackage{graphicx}
\usepackage{amsmath}

\newcommand{\bea}{\begin{eqnarray*}}
\newcommand{\eea}{\end{eqnarray*}}
\newcommand{\be}{\begin{equation}}
\newcommand{\ee}{\end{equation}}
\newcommand{\ben}{\begin{eqnarray*}}
\newcommand{\een}{\end{eqnarray*}}

\date{}
\begin{document}
\title{The minimum number of Hamilton cycles in a hamiltonian threshold graph of a prescribed order\footnote{E-mail addresses:
{\tt 235711gm@sina.com}(P.Qiao),
{\tt zhan@math.ecnu.edu.cn}(X.Zhan). This research  was supported by Science and Technology Commission of Shanghai Municipality (STCSM) grant 13dz2260400 and
 the NSFC grant 11671148.}}
\author{Pu Qiao, Xingzhi Zhan\thanks{Corresponding author.}\\
{\small Department of Mathematics, East China Normal University, Shanghai 200241, China}
 } \maketitle
\begin{abstract}
 We prove that the minimum number of Hamilton cycles in a hamiltonian threshold graph of order $n$ is $2^{\lfloor (n-3)/2\rfloor}$ and this minimum number is
 attained uniquely by the graph with degree sequence $n-1,n-1,n-2,\ldots,\lceil n/2\rceil,\lceil n/2\rceil,\ldots,3,2$ of $n-2$ distinct degrees. This
 graph is also the unique graph of minimum size among all hamiltonian threshold graphs of order $n.$

\end{abstract}


\section{Introduction}

There are few results concerning the precise value of the minimum or maximum number of Hamilton cycles of graphs in a special class with a prescribed order.
For example, it is known that the minimum number of Hamilton cycles in a simple hamiltonian cubic graph of order $n$ is $3,$ which follows from Smith's theorem [1, p.493]
and an easy construction [10, p.479], but the maximum number of Hamilton cycles is not known; even the conjectured upper bound $2^{n/3}$ [2, p.312] has not been proved. Another example is Sheehan's conjecture that every simple hamiltonian 4-regular graph has at least two Hamilton cycles [10] (see also [1, p.494 and p.590]), which is still unsolved.

In this paper we will determine the minimum number of Hamilton cycles in a hamiltonian threshold graph of order $n$ and the unique minimizing graph. Threshold graphs were introduced by
Chv$\acute{a}$tal and Hammer [3] in 1973. Besides the original definition, seven equivalent characterizations are given in the book [7].

{\bf Definition 1.} A finite simple graph $G$ is called a {\it threshold graph} if there exists a nonnegative real-valued function $f$ defined on the vertex set of $G,$
$f: V(G)\rightarrow {\mathbb R}$ and a nonnegative real number $t$ such that for any two distinct vertices $u$ and $v,$ $u$ and $v$ are adjacent if and only if
$f(u)+f(v)>t.$

The class of threshold graphs play a special role for many reasons, some of which are the following: 1) They have geometrical significance. Let $\Omega_n$ be the convex hull of all degree sequences of the simple graphs of order $n.$ Then the extreme points of the polytope $\Omega_n$ are exactly the degree sequences of threshold graphs of order $n$ [6] (for another proof see [9]). 2) A nonnegative integer sequence is graphical if and only if it is majorized by the degree sequence of some threshold graph [9]. 3) A graphical sequence has a unique labeled realization if and only if it is the degree sequence of a threshold graph [7, p.72].

For terminology and notations we follow the textbooks [1,11]. The order of a graph is its number of vertices, and the size its number of edges. We regard isomorphic graphs as the same graph. Thus for two graphs $G$ and $H,$ $G=H$ means that $G$ and $H$ are isomorphic. $N(v)$ and $N[v]$ denote the neighborhood and closed neighborhood of a vertex $v$ respectively. For a real number $r,$ $\lfloor r\rfloor$ denotes the largest integer less than or equal to $r,$ and $\lceil r\rceil$ denotes the least integer larger than or equal to $r.$ The notation $|S|$ denotes the cardinality of a set $S.$

\section{Main Results}

Let $G=(V,E)$ be a graph whose distinct positive vertex-degrees are $\delta_1<\cdots<\delta_m$ and let $\delta_0=0.$ Denote $D_i=\{v\in V| {\rm deg}(v)=\delta_i\}$
for $i=0,1,\ldots,m.$ The sequence $D_0,D_1,\ldots,D_m$ is called the {\it degree partition} of $G.$ Each $D_i$ is called a {\it degree set.} Sometimes when $D_0$ is empty it may be omitted. These notations will be used throughout. We will need the following characterization [7, p.11] which describes the basic structure of a threshold graph.

{\bf Lemma 1.} {\it $G$ is a threshold graph if and only if for each $v\in D_k,$
$$N(v)=\bigcup_{j=1}^k D_{m+1-j}\quad {\rm if }\,\,\, k=1,\ldots,\lfloor m/2\rfloor$$
$$N[v]=\bigcup_{j=1}^k D_{m+1-j}\quad {\rm if }\,\,\, k=\lfloor m/2\rfloor+1,\ldots,m.$$
In other words, for $x\in D_i$ and $y\in D_j,$ $x$ is adjacent to $y$ if and only if $i+j>m.$}

Clearly, Lemma 1 not only implies another characterization that the vicinal preorder of a threshold graph is a total preorder, but also indicates that
every threshold graph is determined uniquely by its degree sequence [7, p.72].

The following lemma can be found in [7, pp.11-13].

{\bf Lemma 2.} {\it For any threshold graph,
$$\delta_{k+1}=\delta_k +|D_{m-k}| \quad {\rm for}\,\,\, k=0,1,\ldots,m,\,\,k\neq \lfloor m/2\rfloor$$
$$\delta_{k+1}=\delta_k +|D_{m-k}|-1 \quad {\rm for}\,\,\, k=\lfloor m/2\rfloor.$$
}

For two subsets $S$ and $T$ of the vertex set of a graph $G,$ the notation  $[S, \,T]$ denotes the set of edges of $G$ with one end-vertex in $S$ and the other end-vertex
in $T.$ Here $S$ and $T$ need not be disjoint. In the case $T=S,$ $[S,\,S]$ is just the edge set of the subgraph  $G[S]$ of $G$  induced by $S.$
Next we define a new concept which will be used in the proofs.

{\bf Definition 2.} An edge of a threshold graph $G$ with degree partition $D_0,D_1,\ldots,D_m$  is called a {\it key edge} of $G$ if it
lies in $[D_k,\,D_{m+1-k}]$ for some $k$ with $1\le k\le \lceil m/2\rceil.$

Thus when $m$ is even we have only one type of key edges, and when $m$ is odd ($m\ge 3$) we have two types of key edges.
For example, if $m=4$ then the set of key edges is $[D_1,\,D_4]\cup [D_2,\,D_3]$ while if $m=5$ then the set of key edges is 
$[D_1,\,D_5]\cup [D_2,\,D_4]\cup [D_3,\,D_3].$
We will need the following two lemmas concerning properties of key edges.

{\bf Lemma 3.} {\it If $e$ is a key edge of a threshold graph $G,$ then $G-e$ is a threshold graph.}

{\bf Proof.} Denote $G^{\prime}=G-e$ and let $m^{\prime}$ be the number of distinct positive vertex-degrees of $G^{\prime}.$
Let $e=xy.$ First suppose that $x\in D_j$ and $y\in D_{m+1-j}$ for some $1\le j\le \lfloor m/2\rfloor.$  We write TPO for the conditions in Lemma 1 (suggesting total preorder). To prove that $G^{\prime}$ is a threshold graph, by Lemma 1 it suffices to show that the degree sets of $G^{\prime}$ satisfy TPO. The structural change of the degree partitions depends on the sizes of the two sets $D_j$ and $D_{m+1-j}.$ We distinguish four cases.

Case 1. $|D_j|=1$ and $|D_{m+1-j}|=1.$ The condition $|D_j|=1$ implies that $j=\lfloor m/2\rfloor$ is possible only if $m$ is odd,
since if $m$ is even then $|D_{m/2}|\ge 2.$ Hence $m-j>j,$ implying that $D_{m-j}$ and $D_j$ are two distinct sets. By Lemma 2,
$$\delta_j=\delta_{j-1}+|D_{m+1-j}|=\delta_{j-1}+1\quad {\rm and}\quad \delta_{m+1-j}=\delta_{m-j}+|D_j|=\delta_{m-j}+1. $$
After deleting $e,$ the two sets $D_j$ and $D_{m+1-j}$ become empty, and they disappear in $G^{\prime}.$ $x$ goes to $D_{j-1}$
and $y$ goes to $D_{m-j}.$  Now $m^{\prime}=m-2$ and the adjacency relations among the vertices of $G^{\prime}$ still satisfy TPO.

Case 2. $|D_j|=1$ and $|D_{m+1-j}|\ge 2.$ As in case 1, $D_{m-j}$ and $D_j$ are two distinct sets. By Lemma 2,  we have
$$\delta_{j}=\delta_{j-1}+|D_{m+1-j}|\ge \delta_{j-1}+2\quad {\rm and} \quad \delta_{m+1-j}=\delta_{m-j}+|D_j|=\delta_{m-j}+1.$$
When deleting $e,$ $x$ stays in $D_j$ and $y$ goes to $D_{m-j}.$ Thus $m^{\prime}=m$ and $G^{\prime}$  satisfies TPO.

Case 3. $|D_j|\ge 2$ and $|D_{m+1-j}|=1.$ We have $\delta_j=\delta_{j-1}+|D_{m+1-j}|=\delta_{j-1}+1.$ When deleting $e,$ $x$ goes to $D_{j-1}.$ If $m$ is even, $j=m/2$ and $|D_j|=2,$ then $\delta_{m+1-j}=\delta_{m/2}+|D_{m/2}|-1=\delta_j+1.$ When deleting $e,$ $y$ goes to $D_j$ and the set $D_{m+1-j}$ disappears. Thus $m^{\prime}=m-1.$ In all other cases, we have $\delta_{m+1-j}\ge \delta_{m-j}+2.$ In fact, if $m$ is odd or $m$ is even and $j<m/2,$ we have
$\delta_{m+1-j}=\delta_{m-j}+|D_j|\ge\delta_{m-j}+2,$ while if $m$ is even, $j=m/2$ and $|D_j|\ge 3,$ we have
$\delta_{m+1-j}=\delta_{m-j}+|D_j|-1\ge\delta_{m-j}+2.$ When deleting $e,$ $y$ remains in $D_{m+1-j}.$ Thus $m^{\prime}=m.$ In each case,
$G^{\prime}$  satisfies TPO.

Case 4. $|D_j|\ge 2$ and $|D_{m+1-j}|\ge 2.$ We have $\delta_{j}=\delta_{j-1}+|D_{m+1-j}|\ge \delta_{j-1}+2.$ If $m$ is even, $j=m/2$ and $|D_j|=2,$ then
$\delta_{m+1-j}=\delta_{j+1}=\delta_j+|D_j|-1=\delta_j+1.$ When deleting $e,$ $x$ remains in $D_j$ (but with degree $\delta_j -1$) and a new degree set
$\{y\}\cup(D_j\setminus\{x\})$ appears. Now $m^{\prime}=m+1.$ In all other cases, two new degree sets appear, one containing only $x$ and the other containing only
$y,$ so that $m^{\prime}=m+2.$ In either case, $G^{\prime}$  satisfies TPO and hence it is a threshold graph.

Now suppose that $m$ is odd and $x,y\in D_t$ where $t={\lfloor m/2\rfloor+1}=\lceil m/2\rceil.$ Apply Lemma 2. If $|D_t|=2,$ when deleting $e,$ both $x$ and $y$ go
to $D_{\lfloor m/2\rfloor}$ and the degree set $D_t$ disappears. Then $m^{\prime}=m-1$ and $G^{\prime}$ satisfies TPO. Otherwise $|D_t|\ge 3.$ When deleting $e,$ a
new degree set $\{x,y\}$ appears, where $x$ and $y$ are nonadjacent. In this case $m^{\prime}=m+1$ and $G^{\prime}$ again satisfies TPO.$\Box$

{\bf Lemma 4.} {\it Every key edge of a hamiltonian threshold graph lies in at least one Hamilton cycle.}

{\bf Proof.} Let $G$ be a hamiltonian threshold graph with degree partition $D_1,\ldots,D_m.$ Let $e=xy$ be a key edge of $G$ with $x\in D_j$ and $y\in D_{m+1-j}$
for some $1\le j\le \lceil m/2\rceil.$ Choose any Hamilton cycle $C$ of $G.$ If $e$ lies in C, we are done. Otherwise let $C=(x,s,\ldots,y,t,\ldots).$ Then
$s$ and $x$ are adjacent, and $t$ and $y$ are adjacent. Applying Lemma 1 we deduce that $s$ and $t$ are adjacent. Now the classical cycle exchange [1, p.485]
with $x^{+}=s$ and $y^{+}=t$ yields a new Hamilton cycle containing the edge $e.$$\Box$

Different necessary and sufficient conditions for a threshold graph to be hamiltonian are given by Golumbic [4], Harary and Peled [5], and Mahadev and Peled [8].
What we need is the following one by Golumbic [4, p.231] whose proof can be found in [7, p.25].

{\bf Lemma 5.} {\it Let $G$ be a threshold graph of order at least $3$ with the degree partition $D_0,D_1,\ldots,D_m.$ Then $G$ is hamiltonian
if and only if $D_0=\phi,$
$$\sum_{j=1}^k |D_j|<\sum_{j=1}^k |D_{m+1-j}|, \quad k=1,\ldots, \lfloor (m-1)/2\rfloor$$
and if $m$ is even, then $\sum_{j=1}^{m/2} |D_j|\le \sum_{j=1}^{m/2} |D_{m+1-j}|.$}

{\bf Definition 3.} For every integer $n\ge 3,$ we denote by $G_n$ the graph with degree sequence $n-1,n-1,n-2,\ldots,\lceil n/2\rceil,\lceil n/2\rceil,\ldots,3,2$ of $n-2$ distinct degrees.

$G_n$ is a hamiltonian threshold graph. $G_8$ is depicted in Figure 1.

\vskip 3mm
\par
 \centerline{\includegraphics[width=0.6\textwidth]{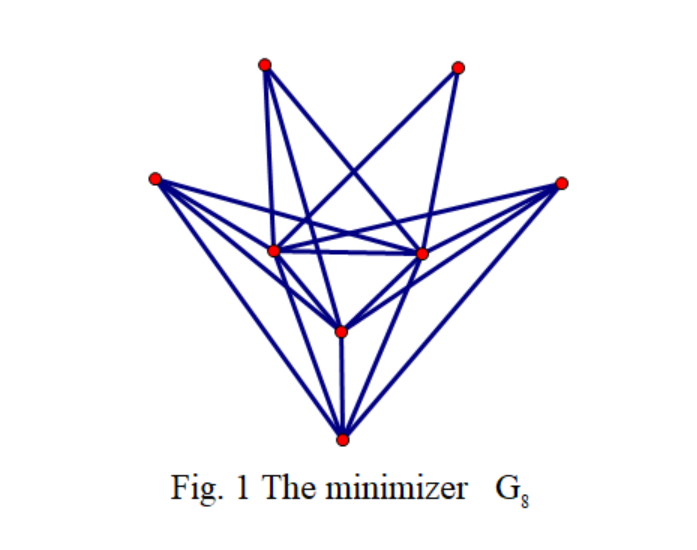}}
\par

Now we are ready to prove the main results.

{\bf Theorem 6.} {\it The minimum number of Hamilton cycles in a hamiltonian threshold graph of order $n$ is $2^{\lfloor (n-3)/2\rfloor}$ and this minimum
number is attained uniquely by the graph $G_n.$}

{\bf Proof.} We first determine the minimizing graph and then count its number of Hamilton cycles. Let $G$ be a hamiltonian threshold graph of order $n$ having the minimum number of Hamilton cycles. Let $D_1,\ldots,D_m$ be the degree partition of $G.$ Note that for any threshold graph with $m\ge 1,$  we have
$|D_{\lceil m/2\rceil}|\ge 2.$ This follows from
$$1\le \delta_{\lfloor m/2\rfloor +1}-\delta_{\lfloor m/2\rfloor}=|D_{\lceil m/2\rceil}|-1$$
by Lemma 2.

The theorem holds trivially for the case $n=3.$ Next suppose $n\ge 4.$ $m=1$ means that $G$ is a complete graph, which is impossible. Thus $m\ge 2.$ We claim that
$|D_m|=2.$ Lemma 5 with $k=1$ implies $|D_m|\ge 2.$ Hence it suffices to prove $|D_m|\le 2.$ To the contrary suppose $|D_m|\ge 3.$ Let $e$ be any edge in $[D_1,\,D_m].$ Then $e$ is a key edge by definition. By Lemma 3, $G-e$ is a threshold graph. Since $G$ is a hamiltonian threshold graph, its degree sets $D_1,\ldots,D_m$ satisfy the inequalities in Lemma 5. Analyzing the change of degree partitions from $G$ to $G-e$ as in the proof of Lemma 3, we see that the degree sets of $G-e$ also satisfy the inequalities in Lemma 5. Hence by Lemma 5,
$G-e$ is hamiltonian. Deleting any edge cannot increase the number of Hamilton cycles. By Lemma 4, the key edge $e$ lies in at least one Hamilton cycle of $G.$
It follows that $G-e$ has fewer Hamilton cycles than $G,$ contradicting the minimum property of $G.$ This proves $|D_m|=2.$

If $m=2,$ then by Lemma 5 we have $|D_1|\le 2.$ Since $n\ge 4,$ we must have $n=4$ and $|D_1|=2.$ Then applying Lemma 1 we deduce that $G$ has the degree sequence $3,3,2,2,$ so that $G=G_4.$ Next suppose $m\ge 3.$ By Lemma 5, $|D_1|<|D_m|=2.$ Hence $|D_1|=1.$ We first consider the case $m\ge 4$ (The case $m=3$ will be treated later). We claim that $|D_{m-1}|=1.$ Otherwise, as argued above, deleting any key edge $f$ in $[D_2,\,D_{m-1}]$ would reduce the number of Hamilton cycles such that $G-f$ is still a hamiltonian threshold graph, a contradiction. Then using the fact that $|D_1|=1$ and $|D_m|=2$ and applying Lemma 5
we deduce that if $m$ is odd or if $m$ is even and $m\ge 6$ then $|D_2|=1,$ and if $m=4$ then $|D_2|=2.$ Continuing in this way, by successively deleting a key edge in $[D_j,\,D_{m+1-j}]$ for $j=2,\dots,\lfloor m/2\rfloor$ if $|D_{m+1-j}|\ge 2$ we conclude that $|D_{m+1-j}|=1$ for each $j=2,\dots,\lfloor m/2\rfloor.$ Then using the fact that $|D_1|=1$ and $|D_m|=2$ and applying Lemma 5, we conclude that $|D_i|=1$ for each $i=2,\dots,\lfloor m/2\rfloor -1$ and that if $m$ is odd then $|D_{\lfloor m/2\rfloor}|=1$ and if $m$ is even then $|D_{m/2}|=2.$ Thus, if $m$ is even then $n=m+2$ is even, $G$ has the degree sequence
$n-1,n-1,\ldots,n/2,n/2,\ldots,3,2$ and hence $G=G_n.$

If $m\ge 3$ and $m$ is odd, we assert that $|D_{\lceil m/2\rceil}|=2.$  As remarked at the beginning, we always have  $|D_{\lceil m/2\rceil}|\ge 2.$ Thus it suffices to show $|D_{\lceil m/2\rceil}|\le 2.$ To the contrary suppose $|D_{\lceil m/2\rceil}|\ge 3.$ By Lemma 4, any key edge $h$ in $G[D_{\lceil m/2\rceil}]$ lies in at least one Hamilton cycle. With the assumption that $|D_{\lceil m/2\rceil}|\ge 3,$  applying Lemma 3 and Lemma 5 we see that $G-h$ is also a hamiltonian threshold graph with fewer Hamilton cycles than $G,$ a contradiction. This shows $|D_{\lceil m/2\rceil}|=2.$ Now $n=m+2$ is odd. Combining all the above information about $G$ we deduce that $G$ has the degree sequence $n-1,n-1,\ldots,(n+1)/2,(n+1)/2,\ldots,3,2$ and hence $G=G_n.$

Denote the number of Hamilton cycles of $G_n$ by $f(n).$ Since $f(3)=f(4)=1,$ to prove $f(n)=2^{\lfloor (n-3)/2\rfloor}$ it suffices to show the following

Claim. For every integer $k\ge 2,$
$$f(2k-1)=f(2k) \quad {\rm and} \quad f(2k+1)=2f(2k). $$

In $G_{2k},$ let $D_k=\{x,y\}$ and $D_{k+1}=\{z\}.$ By Lemma 5, neither $G_{2k}-xz$ nor $G_{2k}-yz$ is hamiltonian.
Thus the path $xzy$ must lie in every Hamilton cycle of $G_{2k}.$ Deleting the vertex $z$ and adding the edge $xy$ we obtain a graph which is isomorphic
to $G_{2k-1}$ and has the same number of Hamilton cycles as $G_{2k}.$ Hence $f(2k-1)=f(2k).$

In $G_{2k+1}$, let $D_{k+1}=\{u,v\}.$ Then the edge $uv$ lies in every Hamilton cycle of $G_{2k+1}.$ Denote $G^{\prime}=G-v.$
Clearly $G^{\prime}$ is isomorphic to $G_{2k}$ and hence $G^{\prime}$ has $f(2k)$ Hamilton cycles. Since $N(u)\setminus \{v\}=N(v)\setminus \{u\}$ in $G,$
from each Hamilton cycle of $G^{\prime}$ we can obtain two distinct Hamilton cycles of $G_{2k+1}$ by replacing the vertex $u$ by the edge $uv$ in two ways.
More precisely, a Hamilton cycle $(\ldots,s,u,t,\ldots)$ of $G^{\prime}$ yields two Hamilton cycles $(\ldots,s,u,v,t,\ldots)$ and $(\ldots,s,v,u,t,\ldots)$
of $G.$ Conversely every Hamilton cycle of $G_{2k+1}$ can be obtained in such a vertex-to-edge expansion from a Hamilton cycle of $G^{\prime}.$ Hence $f(2k+1)=2f(2k).$ This shows the claim and completes the proof.$\Box$

The above proof of Theorem 6 also proves that $G_n$ is the unique graph that has the minimum size among all hamiltonian threshold graphs of order $n.$
To see this, just replace the assumption that $G$ has the minimum number of Hamilton cycles by the one that $G$ has the minimum size. Also note that the size
of a threshold graph is easy to count, since it is a split graph with the clique $\bigcup_{j=\lfloor m/2\rfloor +1}^m D_j$ and the independent set
$\bigcup_{j=1}^{\lfloor m/2\rfloor}D_j.$ Thus we have the following result.

{\bf Theorem 7.} {The minimum size of a hamiltonian threshold graph of order $n$ is
$$\begin{cases} (n^2+2n-3)/4 \quad {\rm if}\,\,\,n\,\,{\rm is}\,\,{\rm odd}\\
         (n^2+2n-4)/4 \quad {\rm if}\,\,\,n\,\,{\rm is}\,\,{\rm even}\end{cases}$$
and this minimum size is attained uniquely by the graph $G_n.$
}


\begin{thebibliography}{99}
\bibitem{1} J.A. Bondy and U.S.R. Murty, Graph Theory, GTM 244, Springer, 2008.
\bibitem{2} G.L. Chia and C. Thomassen, On the number of longest and almost longest cycles in cubic graphs, Ars Combinatoria, 104(2012), 307-320.
\bibitem{3} V. Chv$\acute{a}$tal and P.L. Hammer, Set-packing problems and threshold graphs, CORR 73-21, University of Waterloo, Canada, August 1973.
\bibitem{4} M.C. Golumbic, Algorithmic Graph Theory and Perfect Graphs, Second Edition, Elsevier, 2004.
\bibitem{5} F. Harary and U.N. Peled, Hamiltonian threshold graphs, Discrete Appl. Math., 16(1987), 11-15.
\bibitem{6} M. Koren, Extreme degree sequences of simple graphs, J. Combin. Theory Ser. B, 15(1973), 213-224.
\bibitem{7} N.V.R. Mahadev and U.N. Peled, Threshold Graphs and Related Topics, Elsevier Science B.V., 1995.
\bibitem{8} N.V.R. Mahadev and U.N. Peled, Longest cycles in threshold graphs, Discrete Math., 135(1994), 169-176.
\bibitem{9} U.N. Peled and M.K. Srinivasan, The polytope of degree sequences, Linear Algebra Appl., 114/115(1989), 349-377.
\bibitem{10} J. Sheehan, The multiplicity of hamiltonian circuits in a graph, in Recent Advances in Graph Theory, 477-480, Academia, Prague, 1975.
\bibitem{11} D.B. West, Introduction to Graph Theory, Prentice Hall, Inc., 1996.
\end{thebibliography}
\end{document}